\documentclass[a4paper,twoside]{article}
\usepackage[T1]{fontenc}
\usepackage[reqno,tbtags]{amsmath}
\usepackage{amssymb} 
\usepackage{amsfonts}
\usepackage{amsthm}
\usepackage{graphicx}
\usepackage{fancyhdr}

\newcommand{\N}{\mathbb{N}}
\newcommand{\Z}{\mathbb{Z}}

\newcommand{\R}{\mathbb{R}}

\newcommand{\stable}{\mathrm{s}}
\newcommand{\unstable}{\mathrm{u}}

\newcommand{\floor}[1]{\floor #1 \floor}

\newcommand{\vect}[1]{\overline{#1}}
\newcommand{\diff}[1]{\, \mathrm{d} #1}
\newcommand{\seq}[1]{\underline{#1}}
\newcommand{\sbr}{\mathrm{SBR}}

\theoremstyle{plain}
\newtheorem{theorem}{Theorem}
\newtheorem{lemma}{Lemma}
\newtheorem{corollary}{Corollary}

\newtheorem*{remark}{Remark}

\begin{document}

  \title{Absolutely continuous invariant measures for some piecewise hyperbolic affine maps}
  \author{Tomas Persson\\
{\footnotesize \sl Centre for Mathematical Sciences, Department of
Mathematics,} \\ {\footnotesize \sl Lund Institute of Technology,
P.O. Box 118, SE-22100 Lund, Sweden}\\
{\footnotesize {\tt tomasp@maths.lth.se}, {\tt http://www.maths.lth.se/\~{}tomasp} }}

  \maketitle

  \begin{abstract}
    A class of piecewise affine hyperbolic maps on a bounded subset of the plane is considered.
    It is shown that if a map from this class is sufficiently area-expanding then almost surely
    this map has an absolutely continuous invariant measure.
  \end{abstract}

  \noindent {\bf Acknowledgement.} The author is very grateful to J\"{o}rg Schmeling for interesting
    and useful discussions.

  \section{Introduction}
    In \cite{Pesin}, Pesin studied a general class of piecewise
    diffeomorphisms with a hyperbolic attractor. He showed the existence
    of the Sinai-Bowen-Ruelle measure, or $\sbr$-measure for short, and
    studied the ergodic properties of this measure. If $f: M \rightarrow
    M$ is the system in question then the $\sbr$-measure is a weak limit
    point of the sequence of measures
    \[
      \mu_n = \frac{1}{n} \sum_{k=0}^{n-1} \nu \circ f^{-k},
    \]
    where $\nu$ denotes the Lebesgue measure. Pesin showed that
    the $\sbr$-measure has at most countably many ergodic
    components.  This measure is the
    physically relevant measure as it captures the behaviour of the orbits
    of points from a set of positive Lebesgue measure.

    For a more restricted class, Sataev \cite{Sataev} showed
    that there are only finitely many ergodic components. 
    Schmeling and Troubetzkoy studied in \cite{Schmeling-Troubetzkoy} a more general
    class than Pesin's and proved the existence of the
    $\sbr$-measure. Their method to deal with the non-invertibility of
    the system was to lift the system to a higher dimension and get an
    invertible system on which the calculations were
    made. In this way methods from invertible systems could be used. The
    result could then be projected back to the original system.

    In \cite{Ale-Yor}, Alexander and Yorke considered a one parameter
    class of maps called the fat baker's transformations. These maps are
    piecewise affine maps of the square with one expanding and one
    contracting direction. Their results together with the result of
    Solomyak in \cite{Sol}, imply that for a positive measure set of
    parameters, there is an absolutely 
    continuous invariant measure.

    The Belykh map, was first introduced in \cite{Belykh} by
    Belykh. Schmeling and \mbox{Troubetzkoy} considered in
    \cite{Schmeling-Troubetzkoy} the Belykh map for a wider range of
    parameters. The fat baker's transformations are a
    special case of  the Belykh map in this wider range of parameters.
    The Belykh map was further investigated in \cite{Joerg_belykh}.

    In this article we consider a class of piecewise affine hyperbolic maps
    on a set $K \subset \R^2$, with one contracting and one expanding direction.
    This class is contained in the class of maps studied in \cite{Schmeling-Troubetzkoy}
    and it contains the Belykh maps as well as the fat baker's transformations.

    It is shown that if a functions from this class is sufficiently area-expanding then 
    almost surely (in the sense of Corollary \ref{cor:acim}) there is an absolutely continuous 
    invariant measure. The method used to show this is a development of the method from 
    \cite{Per-Sol}. Here a new problem arises: The symbolic space changes as the parameters changes.
    In this paper a way to handle this problem is introduced. The different symbolic spaces are
    embedded in a larger space and certain estimates are carried out that makes this larger space
    possible to handle.

    In \cite{Tsujii_fat} and \cite{Tsujii_hyp}, Tsujii considered two classes of maps in two dimension 
    and showed that almost all of these maps have absolutely continuous invariant measure. These two
    classes are different from the class of maps considered in this paper. Tsujii also used the method
    from \cite{Per-Sol}, but in a different way than is used in this paper. 

    Similar results in two dimensions, but in the case of expanding
    maps, were independently obtained by Buzzi in \cite{Buzzi_2dim} and
    Tsujii in \cite{Tsujii_2dim}. The corresponding results for arbitrary
    dimension are in \cite{Buzzi_highdim} and \cite{Tsujii_highdim}.

  \section{A class of piecewise hyperbolic maps} \label{sec:class}
    Let $K \subset \R^2$ be compact and connected. Assume that $K$ can be decomposed according to
    \[
      K = \bigcup_{i=1}^a \overline{K}_i
    \]
    where each $K_i$ is an open and non-empty set with the boundary consisting of finitely many $C^2$ curves.
    Thus, there are closed $C^2$ curves $N_{i}$ and $M_{i}$ such that
    \[
      \bigcup_{i=1}^a \partial K_i = \Bigl( \bigcup_{i=1}^{b} N_i \Bigr) \cap \Bigl( \bigcup_{i=1}^c M_i
      \Bigr), \
      \partial K = \bigcup_{i=1}^c M_i.
    \]
    Let $\mathcal{Z} = \{K_i\}$ denote the partion of $K$.

    Assume that the sets $N_i \cap N_j$, $M_i \cap M_j$ and $N_i \cap M_k$ consists of finitely many points if
    $i \ne j$,
    and there exists a constant $H$ such that if $(t_1,t_2) \in T_p N_{ij}$ then $|t_2/t_1| < H$. 
    Let $N = \cup N_i$ and $M = \cup M_i$.

\begin{figure}[!ht]
      \begin{center}
        \includegraphics[scale=0.8]{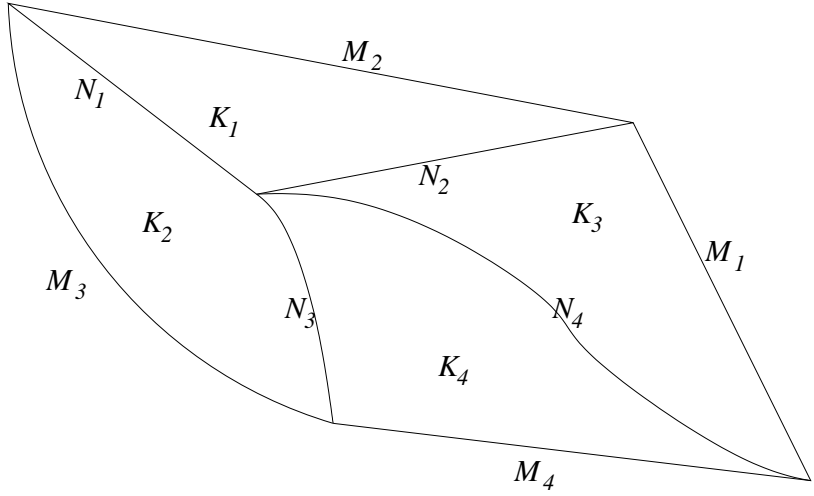}
      \end{center}
      \caption{An example of the domain $K$.}
      \label{fig:domain}
    \end{figure}

    See Figure \ref{fig:domain} for an example of $K$, $N_i$
    and $M_i$.

    Consider maps $f: K \setminus N \rightarrow K$ that satisfy the following two conditions, (A1)
    and (A2).

    \begin{itemize}
      \item[(A1)]  There are numbers 
      $\lambda_1, \ldots, \lambda_a < 1 < \gamma_1, \ldots, \gamma_a$ and $u_1, \ldots u_a, v_1,
      \ldots v_a \in \R$ with $u_i \ne u_j$ whenever $i \ne j$,
      such that for any $i = 1,\ldots, a$ the map $f$ restricted to $K_i$ is defined by
      \[
        f|_{K_i} (x_1,x_2) = f_i (x_1,x_2) = (\lambda_i x_1 + u_i, \gamma_i x_2 + v_i).
      \]
    \end{itemize}

    The notation $f_{\vect{\lambda}, \vect{\gamma}, \vect{u}, \vect{v}}$ for $f$ will be used to
    emphasise 
    the dependence on the parameters $\vect{\lambda} = (\lambda_1, \ldots, \lambda_a)$, 
    $\vect{\gamma} = (\gamma_1, \ldots, \gamma_a)$, 
    $\vect{u} = (u_1, \ldots u_a)$ and $\vect{v} = (v_1, \ldots v_a)$.

    Let $\mathcal{Z}_0 = \mathcal{Z}$ and assume that $\mathcal{Z}_k = \{K_i^{(k)}\}_{i=1}^{a_k}$ is defined.
    (Note that $a_0 = a$.)
    Then define the partition $\mathcal{Z}_{k+1} = \{K_i^{(k+1)}\}_{i=1}^{a_{k+1}}$ by
    \[
      \mathcal{Z}_{k+1} = \mathcal{Z}_k \wedge \{ f (K_i^{(k)})\}_{i=1}^{a_{k}}.
    \]
    The set $\bigcup_{i=1}^{a_k} K_i^{(a_k)}$ is the set of points $x \in K$ such that for each
    $l = 0,1,\ldots,k$ the point $f^l (x)$ is defined and $f^l (x) \not \in N \cup M$.

    Since each $K_i$ has piecewise $C^2$ boundary so has each $K_j^{(k)}$. There are thus closed
    $C^2$ curves $N_i^{(k)}$ such that
    \[
      \bigcup_{i=1}^{a_k} \partial K_i^{(k)} = \bigcup_{i=1}^{b_k} N_i^{(k)}
    \]
    and $N_i^{(k)} \cap N_j^{(k)}$ is a finite set if $i \ne j$.

    \begin{itemize}
      \item[(A2)]
        There is a number $\tau \geq 1$ such that $(\min \gamma_i)^\tau > D_\tau + 1$ where
        \[
          D_\tau = \max \{ \# A \mid A \subseteq \{1,2, \ldots, b_\tau \} \text{ s.t. } 
          \bigcap_{i \in A} N_i^{(\tau)} \ne \emptyset \}.
        \]
        That is, $D_\tau$ is the maximal number of lines from $\{N^{(\tau)}_i\}$ that crosses at
        one point.
        The number $D_\tau$ is finite since the set $\{N_i^{(\tau)} \}$ is finite.   
    \end{itemize}

    \begin{remark}
      Condition (A2) implies that the multiplicity entropy (see \cite{Buzzi_intrinsic} and \cite{Krug-Ryp}
      for a definition) is less than the positive Lyapunov exponent.
    \end{remark}

  \section{The results}

    We are going to prove the following theorem.

    \begin{theorem} \label{the:acim}
      Assume that for all $t \in I =(t_0,t_1)$, the maps 
      $f_{t \vect{\lambda}, \vect{\gamma}, \vect{u}, \vect{u}}$ satisfy the conditions (A1) and (A2)
      with a uniform $\tau$. If
      \[
        \frac{t_0 \min \{\lambda_i\} \min \{\gamma_i^2\}}{\max \{\gamma_i\}} > 1,
      \]
      and one of the following conditions are satisfied
      \begin{align*}
        &\text{1.}& & t_1 \lambda_{\max} < 0.5,& & \frac{\max \{|u_i - u_j|, |u_i|\}}
          {\min \{|u_i -u_j|: u_i \ne u_j\}} <
          \frac{1 - t_1 \lambda_{\max}}{t_1 \lambda_{\max} - 2 (t_1 \lambda_{\max})^3}, \\
        &\text{2.}& & t_1 \lambda_{\max} < 0.61,& & \frac{\max \{|u_i - u_j|, |u_i|\}}
          {\min \{|u_i -u_j|: u_i \ne u_j\}} <
          \frac{1 - t_1 \lambda_{\max}}{t_1 \lambda_{\max}- 2 ( t_1 \lambda_{\max})^4}, \\
        &\text{3.}& & t_1 \lambda_{\max} < 0.68, & & \frac{\max \{|u_i - u_j|, |u_i|\}}
          {\min \{|u_i -u_j|: u_i \ne u_j\}} <
          \frac{1 - t_1 \lambda_{\max}}{t_1 \lambda_ {\max}- 2 ( t_1 \lambda_{\max})^5}, 
      \end{align*}
      where $\lambda_{\max} = \max\{\lambda_i\}$, then for almost every $t\in I$ there exists an 
      $f_{t \vect{\lambda}, \vect{\gamma}, \vect{u}, \vect{v}}$-invariant measure,
      absolutely continuous with respect to Lebesgue measure.
    \end{theorem}

    Conclude the following

    \begin{corollary} \label{cor:acim}
      Let $P$ be the set of parameters $(\vect{\lambda}, \vect{\gamma}, \vect{u}, \vect{v})$ such that
      \[
        \frac{ \min \{\lambda_i\} \min \{\gamma_i^2\}}{\max \{\gamma_i\}} > 1,
      \]
      and one of the following conditions are satisfied
      \begin{align*}
        &\text{1.}& &\max \{\lambda_i\} < 0.5,& & \frac{\max \{|u_i - u_j|, |u_i|\}}
          {\min \{|u_i -u_j|: u_i \ne u_j \}} <
          \frac{1 - \max \{\lambda_i\}}{\max\{\lambda_i\} - 2 (\max \{\lambda_i\})^3}, \\
        &\text{2.}& &\max \{\lambda_i\} < 0.61,& & \frac{\max \{|u_i - u_j|, |u_i|\}}
          {\min \{|u_i -u_j|: u_i \ne u_j\}} <
          \frac{1 - \max \{\lambda_i\}}{\max\{\lambda_i\} - 2 (\max \{\lambda_i\})^4}, \\
        &\text{3.}& &\max \{\lambda_i\} < 0.68,& & \frac{\max \{|u_i - u_j|, |u_i|\}}
          {\min \{|u_i -u_j|: u_i \ne u_j\}} <
          \frac{1 - \max \{\lambda_i\}}{\max\{\lambda_i\} - 2 (\max \{\lambda_i\})^5}. 
      \end{align*} 
      If $f_{\vect{\lambda}, \vect{\gamma}, \vect{u}, \vect{u}}$ satisfies the conditions (A1) and (A2)
      then for Lebesgue almost every $(\vect{\lambda}, \vect{\gamma}, \vect{u}, \vect{v}) \in P$,
      there is an $f_{\vect{\lambda}, \vect{\gamma}, \vect{u}, \vect{v}}$-invariant 
      measure, absolutely continuous with respect to Lebesgue measure.    
    \end{corollary}

  \section{A condition on transversality for power series}
    The corollary of the following lemma will be used to prove Theorem \ref{the:acim}. The Lemma appears 
    in a somewhat less general form in \cite{Per-Sol}. The proof from \cite{Per-Sol} works here as well.

    \begin{lemma} \label{lem:trans}
      Let $C \geq 1$. Then there is a constant $\delta > 0$ such that for any function $g$
      of the form
      \begin{equation} \label{g}
        g (x) = 1 + \sum_{k=1}^\infty b_k x^k, \ b_k \in [-C,C],
      \end{equation}
      the following implication holds true for $n=2,3,4$
      \[
        g(x) \leq \delta, \ x \in (0, Q_n),\ C < f_n(x)
        \Longrightarrow g' (x) \leq -\delta,
      \]
      where $f_n (x) = \frac{1 - x}{x - 2 x^{n+1}}$, $Q_2 = 0.5$, $Q_3 = 0.61$ and $Q_4 = 0.68$.
    \end{lemma}

    \begin{proof}
      Let $h_n (x) = 1 - C \sum_{k=1}^n x^k + C \sum_{k=n+1}^\infty x^k = 1 - C \frac{x - 2 x^{n+1}}{1-x}$ and
      $f_n(x) = \frac{1 - x}{x - 2 x^{n+1}}$.
      Let $Q_2 = 0.5$, $Q_3 = 0.61$ and $Q_4 = 0.68$.
      This implies that there is a $\delta > 0$ such that for $n=2,3,4$
      \begin{itemize}
        \item[i)] $h_n (x) > \delta$ if $C < f_n (x)$
        \item[ii)] $h_n'(x) < - \delta$ if $x \in (0, Q_n)$.
      \end{itemize}
      Thus $h_n(x) > \delta$ and $h_n'(x) < - \delta$ provided $x < Q_n$ and $C < f_n (x)$.

      Let $g(x)$ be of the form \eqref{g}. There are $c_k \geq 0$ such that the function
      $G_n(x) = g(x) - h_n(x)$ can be written as $G_n(x) = \sum_{k=1}^n c_k x^k - \sum_{k=n+1}^\infty c_k x^k$.
      By the above mentioned properties of $h_n$ we have for any $x \in (0, Q_n)$, $C < f_n(x)$
      \[
        g(x) < \delta \ \Rightarrow \ G_n(x) < 0 \ \Rightarrow \ G_n'(x) < 0 \ \Rightarrow \ g'(x) < -\delta.
      \]
      The second implication is proved by
      \begin{multline*}
        G_n(x) < 0 \ \Rightarrow \ \sum_{k=1}^n c_k x^k < \sum_{k=n+1}^\infty c_k x^k \ \Rightarrow \\
        \sum_{k=1}^\infty k c_k x^k < \sum_{k=n+1}^\infty k c_k x^k \ \Rightarrow \ G_n'(x) < 0. \qedhere
      \end{multline*}
    \end{proof}

    \begin{corollary} \label{cor:trans}
      Let $n \in \{2, 3, 4\}$ and let $Q_n$ and $f_n$ be as in Lemma \ref{lem:trans}.
      Let $s_k$ be a sequence with $s_k \in [-C,C]$. Then for any $l > 0$ the set
      \[
        \Bigl\{ q \in (q_0, Q_n) \mid C < f_n(q),\
        \Bigl| q^l + \sum_{k=l+1}^\infty s_k q^k \Bigr| < r \Bigr\} 
      \]
      is contained in an interval of length at most $2 \delta^{-1} q_0^{-l} r$.
    \end{corollary}

    \begin{proof}
      Note that $\Bigl| q^l + \sum_{k=l+1}^\infty s_k q^k \Bigr| < r$ implies
      $\Bigl| 1 + \sum_{k=l+1}^\infty s_k q^{k-l} \Bigr| < r q_0^{-l}$ if $q \geq q_0$.
      Lemma \ref{lem:trans} implies that on the
      set $\{q \in (0, Q_n) \mid C < f_n(q)\}$
      the graph of the function $q \mapsto 1 + \sum_{k=l+1}^\infty s_k q^{k-l}$ crosses $0$ transversally
      in at most one point and the slope is at most $- \delta$ around this point.
      Hence $\Bigl| 1 + \sum_{k=l+1}^\infty s_k q^{k-l} \Bigr| < r q_0^{-l}$ on an interval of length not
      more than $2 \delta^{-1} q_0^{-l} r$ and therefore $\Bigl| q^l + \sum_{k=l+1}^\infty s_k q^k \Bigr| < r$
      can only hold in this interval.
    \end{proof}

  \section{Proof of Theorem \ref{the:acim}}
    We will use the method from \cite{Joerg_belykh} to prove Theorem \ref{the:acim}. This method is
    based on \cite{Per-Sol}. The idea is to integrate the density of the measure and then integrate
    with respect to the parameter. If this integral is finite then almost surely the density is integrable
    and so the measure is absolutely continuous with respect to Lebesgue measure. To prove that this is the case
    it is necessary to control how the measure changes with the parameter.

    We begin with some notation and general theory and then make the estimates needed later in the proof.

    \subsection{Notation and general theory} \label{subsec:notation}
    Fix $\vect{\lambda}$, $\vect{\gamma}$, $\vect{u}$ and $\vect{v}$.
    Let $\lambda_{\min} = \min \{\lambda_i\}$, $\lambda_{\max} = \max \{\lambda_i\}$,
    $\gamma_{\min} = \min \{\gamma_i\}$ and $\gamma_{\max} = \max \{\gamma_i\}$.
    We will use the shorter notation $f_t$ to denote
    $f_{t \vect{\lambda}, \vect{\gamma}, \vect{u}, \vect{u}}$.

    Let $\hat{K} = K \times [0,1]$ and $\hat{K}_i = K_i \times [0,1]$. The sets
    $\hat{N}$, $\hat{M},\ldots$ are defined analogously.
    We use the idea from \cite{Schmeling-Troubetzkoy} and lift the map $f_t$ to an
    injective map $\hat{f}_t$ on $\hat{K}$ by
    \[
      \hat{f}_t |_{\hat{K}_i} (x_1,x_2,x_3) = (f_t(x_1,x_2), \theta x_3 + i/(a+1)),
    \]
    where $0 < \theta < 1/(a+1)$. The map $\pi: \hat{K} \rightarrow K, (x_1,x_2,x_3) \mapsto (x_1,x_2)$
    is the projection of
    $\hat{K}$ on $K$. It satisfies $\pi (\hat{f}_t (x_1,x_2,x_3)) = f_t (\pi (x_1,x_2,x_3))$.

    Let $\hat{D}_t^+ = \{\hat{p} \in \hat{K} \mid \hat{f}_t^n (\hat{p}) \not \in \hat{N} \cup \hat{M}, \
    \forall n \in \N \}$ and $\hat{D}_t = \bigcap_{n=0}^\infty \hat{f}_t^n (\hat{D}_t)$.
    The set $\hat{\Lambda}_t = \overline{\hat{D}_t}$ is the attractor of $\hat{f}_t$.

    The condition (A2) is a more general version of condition (H9) in \cite{Pesin}. It appears in
    \cite{Schmeling-Troubetzkoy}. Theorem 6.1 in \cite{Schmeling-Troubetzkoy} can be applied to conclude
    that there are constants $c_t > 0$ such that for any $\varepsilon > 0$ and any $n \in \N$
    \[
      \hat{\nu} ( \hat{f}_t^{-n} (U(\varepsilon, \hat{N} \cup \hat{M}))) < c_t \varepsilon,
    \]
    where $U(\varepsilon, \hat{N} \cup \hat{M})$ denotes the $\varepsilon$-neighborhood of
    $\hat{N} \cup \hat{M}$.
    This shows that the functions $\hat{f}_t$ are in the class of functions from \cite{Pesin}.
    Moreover, the maps
    $f_t$ also satisfies the conditions in \cite{Schmeling-Troubetzkoy}. This gives us the following
    results.

    \begin{itemize}
      \item[i)]
        Let $V^\unstable \subset K$ be a curve in the unstable direction, i.e. there are numbers $\rho, \sigma_1$
        and $\sigma_2$ such that $V^\unstable = \{(x_1, x_2) \in K \mid x_1 = \rho,\ \sigma_1 < x_2 < \sigma_2\}$.
        Let $\hat{V}^\unstable = \pi^{-1} (V^\unstable)$ be the corresponding manifold in $\hat{K}$.
        Let $\nu_{V^\unstable}$ and $\hat{\nu}_{\hat{V}^\unstable}$ denote the normalised Lebesgue measure on
        $V^\unstable$ and $\hat{V}^\unstable$ respectively.
        The sequence of measures $\hat{\mu}_n^t = \frac{1}{n} \sum_{k=0}^{n-1} \hat{\nu}_{\hat{V}^\unstable} \circ 
        \hat{f}_t^{-k}$ 
        converges weakly to an $\sbr$-measure $\hat{\mu}_\sbr^t$. The projection of this measure 
        $\hat{\mu}_\sbr^t \circ \pi^{-1}$ is an $\sbr$-measure for $f_t$ and we thus write 
        $\mu_{\sbr}^t = \hat{\mu}_\sbr^t \circ \pi^{-1}$.

      \item[ii)]
        Given $\varepsilon > 0$, the set
        \[
          \hat{D}_{t, \varepsilon, c} = \{ \hat{x} \in \hat{\Lambda}_t \mid d(\hat{f}_t^n (\hat{x}),
          \hat{N})
          \geq c e^{-\varepsilon n} \}
        \]
        is non-empty if $c$ is sufficiently small and the set $\hat{D}_{t, \varepsilon} = 
        \bigcup_{i=1}^\infty \hat{D}_{t, \varepsilon, i^{-1}}$ has full $\hat{\mu}_\sbr^t$-measure, 
        $\hat{\mu}_\sbr^t (\hat{D}_{t, \varepsilon}) = 1$.  

      \item[iii)]
        The conditional measures of $\hat{\mu}_\sbr^t$ on the unstable manifolds are absolutely
        continuous 
        with respect to Lebesgue measure.

      \item[iv)]
        The entropy of the measure $\hat{\mu}_\sbr^t$ is
        \[
          h_{\hat{\mu}_\sbr^t} = \int_{\hat{\Lambda}_t} \log \chi_t (\hat{x})
          \diff{\hat{\mu}_\sbr^t} (\hat{x}),
        \]
        where $\chi_t (\hat{x})$ is the positive Lyapunov exponent at the point $\hat{x}$ for
        the map $\hat{f}_t$. 
        In particular, $\log (\gamma_{\min}) \leq h_{\hat{\mu}_\sbr^t} \leq \log (\gamma_{\max})$.

      \item[v)] The measure $\hat{\mu}_\sbr^t$ has at most countably many ergodic components.
    \end{itemize}

    The results i) -- v) make it possible to define stable manifolds for 
    $\hat{\mu}_\sbr^t$-almost every $\hat{x} \in \hat{K}$.
    If $\hat{x} = (x_1, x_2, x_3) \in \hat{D}_{t, -\log (t \lambda_{\max})}$ then there is a 
    $c = c(\hat{x})$ 
    such that $\hat{x} \in \hat{D}_{t, -\log (t \lambda_{\max}), c}$, that is 
    $d(\hat{f}_t^n (\hat{x}), \hat{N})
    \geq c (t \lambda_{\max})^n$ for all $n \geq 0$. If $\hat{y} = (y_1, y_2, y_3) \in \hat{K}$ with
    $|y_1 - x_1| < c$ and $y_2 = x_2$ then
    \[
      |\hat{f}_t^n (\hat{y}) - \hat{f}_t^n (\hat{x})|_1 \leq |y_1 - x_1| (t \lambda_{\max})^n 
      < c (t \lambda_{\max})^n \leq d(\hat{f}_t^n (\hat{x}), \hat{N}),
    \]
    where $| \cdots |_1$ denotes the modulus of the difference in the first coordinate. Hence the points 
    $\hat{f}_t^n (\hat{y})$
    and $\hat{f}_t^n (\hat{x})$ are never separated by a discontinuity and we say $\hat{y}$ is in the 
    stable manifold of $\hat{x}$. The stable manifold of $\hat{x}$ is thus defined to be the set
    \[
      \hat{W}^{t, \stable} (\hat{x}) = \{ \hat{y} \in \hat{K} \mid |x_1 - y_1| < c, \ y_2=x_2\}, 
    \]
    where $c$ is the largest constant such that $\hat{x} \in \hat{D}_{t, - \log (t \lambda_{\max}), c}$.
    This defines the stable manifold of $\hat{\mu}_\sbr^t$-a.e. point $\hat{x} \in \hat{K}$ since
    $\hat{\mu}_\sbr^t (\hat{D}_{t,\log (t \lambda_{\max})}) = 1$.

    The stable manifold $W^{t,\stable} (x)$ is defined as the projection of corresponding
    stable manifold in $\hat{K}$.
    All the stable manifolds will therefore be parallel line segments but the length of the manifolds are only
    measurable (\cite{Pesin}).

    Similarly the unstable manifolds can be defined. They consist of
    parallel line segments, orthogonal to the stable manifolds, and their length is measurable.

    The partition of $K$ into stable manifolds is thus measurable
    and the conditional measures on these manifolds  can be defined.
    Take $x \in \Lambda = \pi (\hat{\Lambda})$ and let $\mu_\sbr^{t,\stable,x}$ 
    denote the conditional measure on the stable manifold $W^{t,\stable} (x)$.

    The sequence of measures $\hat{\mu}_n^t$ may converge to a measure which is not ergodic, but if
    $\hat{V}^\unstable \subseteq \hat{W}^{t, \unstable} (\hat{x})$ for some $\hat{x}$ then $\hat{\mu}_n^t$
    will converge to a unique ergodic component. To get control over $\hat{\mu}_\sbr^t$ for 
    almost all $t \in I$ argue as follows.
    Since the size of the unstable manifolds depends only measurable on $t$ there is no $\hat{V}^\unstable$
    that is contained in an unstable manifold for every $t \in I$.
    However, given $\varepsilon_0 > 0$ there is a set $J_0 \subset I$
    with $\nu (I \setminus J_0) < \varepsilon_0$
    and a set $\hat{V}^\unstable$, independent of $t$,
    such that $\hat{V}^\unstable \subseteq \hat{W}^{t, \unstable} (\hat{x})$
    for some $\hat{x}$ whenever $t \in J_0$. This construction has the following consequence.
    For any $t \in J_0$ the measure $\hat{\mu}_n^t$ converges to an ergodic measure and it is not necessary
    to take a subsequence. Indeed, if it is necessary to take a subsequence then there exists a set $A$ such that
    \[
      \liminf_{n\rightarrow \infty} \hat{\mu}_n^t (A) < \limsup_{n\rightarrow \infty} \hat{\mu}_n^t (A).
    \]
    But this would contradict that $\hat{\mu}_n^t$ converges to a unique ergodic component.

    Given a sequence $\{i_n\} \in \{1,2, \ldots, a\}^\Z$ and integers $l, m$ we define the cylinder
    \[
      \hat{C}_{l}^m (\{i_n\},t) = \bigcap_{k=l}^{m} \hat{f}^{-k} (\hat{K}_{i_k})
    \]
    or
    \[
      \hat{C}_{l}^m (i_0,i_1, \ldots, i_{m-l},t) = \bigcap_{k=l}^{m} \hat{f}^{-k} (\hat{K}_{i_{k-l}}).
    \]
    Let $\Sigma_t = \{ \{i_k\} \mid \hat{C}_l^m (\{i_k\},t) \ne \emptyset \ \forall l,m \in \Z\}$ and
    let $\rho_t : \Sigma_t \rightarrow \hat{\Lambda}_t$ be the natural identification of sequences 
    in $\Sigma_t$ and points in $\hat{\Lambda}_t$.

  \subsection{The measure's dependence on the parameters}
    The first step is to estimate how the measure $\hat{\mu}_\sbr^t$ changes with the parameter $t$. This is done by
    using that $\hat{\mu}_{n}^t$ converges weakly to $\hat{\mu}_\sbr^t$.

    Let $L_1 > 0$. The measure $\hat{\mu}_n^t = \frac{1}{n} \sum_{k=0}^{n-1} \hat{\nu}_{\hat{V}^\unstable}
    \circ \hat{f}_t^{-k}$
    converges weakly to $\hat{\mu}_\sbr^t$. Hence there is a number $n_0 (t, L_1)$ such that for any
    cylinder of length $L_1$
    \[
      \frac{1}{2} \leq \frac{\hat{\mu}_\sbr^t (\hat{C}_{-L_1}^0) }{\hat{\mu}_n^t (\hat{C}_{-L_1}^0)} \leq 2,
    \]
    for all $n \geq n_0$. Since $n_0 (t,L_1)$ is measurable with respect to $t$, Lusin's theorem implies that for
    any $\varepsilon_1 > 0$ there is a set $J_1 \subset J_0$ with $\nu (J_0 \setminus J_1) < \varepsilon_1$
    and a number $n_1$ such that $n_1 > n_0 (t, L_1)$ for all $t \in J_1$.

    For a fixed cylinder $\hat{C}_{-L_1}^0$, the measure $\hat{\mu}_{n_0}^t (\hat{C}_{-L_1}^0)$ depends continuously
    on $t$, because the measure $\hat{\mu}_{n_0}^t$ involves taking finitely many preimages with respect to
    $\hat{f}_t$
    and these preimages depends continuously on $t$. It is therefore possible to partition $I$ into
    finitely many subintervals $I = \bigcup_{k = 1}^m I_k$ such that when $t_1, t_2 \in I_k$ then
    \[
      \frac{1}{2} \leq \frac{\hat{\mu}_{n_0}^{t_1} (\hat{C}_{-L_1}^0) }{\hat{\mu}_{n_0}^{t_2} 
      (\hat{C}_{-L_1}^0)} \leq 2,
    \]
    for any cylinder of length $L_1$. This implies that
    \[
      \frac{1}{4} \leq \frac{\hat{\mu}_\sbr^{t_1} (\hat{C}_{-L_1}^0) }{\hat{\mu}_\sbr^{t_2}
      (\hat{C}_{-L_1}^0)} \leq 4,
    \]
    for any cylinder of length $L_1$, provided $t_1, t_2 \in I_k \cap J_1$ for some $k$.

    For each $I_k$, choose a $t_k \in I_k \cap J_1$ and a cylinder $\hat{C}_{-L_1}^0 (\{x_n (t_k)\}, t_k)$ with
    $m_k = \hat{\mu}_\sbr^{t_k} (\hat{C}_{-L_1}^0 (\{x_n (t_k)\}, t_k) > 0$. For every $t \in I_k \cap J_1$, define
    \[
      \hat{\Omega}_{0,t} = \hat{C}_{-L_1}^0 (\{x_n (t_k)\}, t).
    \]
    Then 
    \begin{equation} \label{Omega0}
      \hat{\mu}_\sbr^{t} (\hat{\Omega}_{0,t}) \geq \frac{1}{4} m_k
    \end{equation}
    for all $t \in I_k \cap J_1$.

  \subsection{Entropy}
    It will be necessary to control the number of cylinders and the measure of the cylinders. As already noted,
    the general theory gives that the entropy of the measure $\hat{\mu}_\sbr^t$, satisfy
    $\log (\gamma_{\min}) \leq h_{\hat{\mu}_\sbr^t} \leq \log (\gamma_{\max})$.

    The Shannon-McMillan-Breiman Theorem implies that for $t \in J_1 \cap I_k$ and $\varepsilon_2 > 0$ 
    there is a constant $A(t)$ such that
    \begin{multline*}
      \hat{\mu}_\sbr^t \biggl( \bigcap_{L>0} \bigl\{ \hat{x} \mid 
      \exists \hat{C}_{-L}^0 (\{x_k\},t) \ni \hat{x}
      \text{ with } \\ A(t) (\gamma_{\max} + \varepsilon_2)^{-L} < \hat{\mu}_\sbr^t (\hat{C}_{-L}^0) 
      < A(t) (\gamma_{\min} - \varepsilon_2)^{-L} \bigr\} \biggr) > 1 - \frac{1}{8} m_k.
    \end{multline*}
    The constant $A(t)$ can be chosen so that almost surely this estimate is also valid
    for the conditional measures on the stable manifolds,
    \begin{multline} \label{smb}
      \hat{\mu}_\sbr^t \biggl( \bigcap_{L>0} \bigl\{ \hat{x} \mid 
      \exists \hat{C}_{-L}^0 (\{x_k\},t) \ni \hat{x}
      \text{ with } \\ A(t) (\gamma_{\max} + \varepsilon_2)^{-L} < \hat{\mu}_\sbr^{t, \stable, x} (\hat{C}_{-L}^0) 
      < A(t) (\gamma_{\min} - \varepsilon_2)^{-L} \bigr\} \biggr) > 1 - \frac{1}{8} m_k.
    \end{multline}
    Let $\hat{\Omega}_{\text{SMB}, t}$ be the set whose measure is estimated above.

    An application of Lusin's theorem shows that that given $\varepsilon_3 > 0$ there exists a set
    $J_2 \subset J_1$ and numbers $A_k$ such that
    \begin{align*}
      & \nu (J_2 \setminus J_1) < \varepsilon_3, \\
      & A (t) \leq A_k, \text{ whenever } t \in I_k \cap J_2.
    \end{align*}

    For all $t \in I_k \cap J_2$, define $\hat{\Omega}_t = \hat{\Omega}_{0,t} \cap \hat{\Omega}_{\text{SMB}, t}$.
    The estimates \eqref{Omega0} and \eqref{smb} shows that 
    $\hat{\mu}_\sbr^{t} (\hat{\Omega}_t) \geq \frac{1}{8} m_k$.

    It follows by \eqref{smb} that the number of cylinders of length $L$ in $\hat{\Omega}_t$ satisfy
    \begin{equation} \label{words}
      N_t (n) \leq A_k (\gamma_{\max} + \varepsilon_2)^n.
    \end{equation}

    Let $L_2 > 0$ and consider for each $t \in I_k$ the set of words of length $L_2 + 1$
    \[
      \mathcal{C}_t = \{x_0, x_1, \ldots, x_{L_2} \mid x_i \in \{1,2, \ldots,a\},\ \hat{C}_0^{L_2} (x_0, x_1,
      \ldots, x_{L_2} , t)
      \ne \emptyset \}.
    \]
    The cylinders $\hat{C}_0^{L_2} (\{x_k\}, t)$ change continuously with $t$.
    The fact that the sets $\hat{K}_i$ have piecewise $C^2$ boundaries
    allow us to draw the following conclusion. There is a partition of $I_k$ into
    finitely many intervals $I_{k,l} (L_2)$
    such that $\mathcal{C}_t = \mathcal{C}_{t'}$ if $t', t \in I_{k,l} (L_2)$ for some $l$.
    Indeed, when $t$ runs over $I_k$, each cylinder appears and disappears only finitely many times.
    Let the finite set of $t$, for
    which some cylinder $\hat{C}_0^{L_2} (\{x_k\}, t)$ appears and disappears,
    define the endpoints of the intervals $I_{k,l} (L_2)$.

    Any sequence in $\Sigma_t \cap \rho_t^{-1} (\hat{\Omega}_t)$ can be written as a concatenation of words from 
    $\mathcal{C}_t$.
    Together with \eqref{words} this implies that for each $I_{k,l}$
    the number of words of length $n$ in $\bigcup_{t \in I_{k,l} \cap J_2} \Sigma_t \cap 
    \rho_t^{-1} (\hat{\Omega}_t)$
    does not exceed $\bigl(A_k^\frac{1}{L_2} (\gamma_{\max} + \varepsilon_2) \bigr)^{n+L_2}$ for any $n$.
    Hence, we have the following lemma.

    \begin{lemma} \label{lem:embed}
      For any $\varepsilon_4 > 0$ there is a number $L_3 = L_3 (\varepsilon_4)$ such that if $L_2 > L_3$ then
      for each $I_{k,l} (L_2)$ the symbolic space $\Sigma_{I_{k,l} (L_2)} = \bigcup_{t \in I_{k,l} (L_2) \cap J_2}
      \Sigma_t \cap \rho_t^{-1} (\hat{\Omega}_t)$ satisfy
      \[
        N_{I_{k,l} (L_2)} (n) \leq B_{k,l} (\gamma_{\max} + \varepsilon_2 + \varepsilon_4)^n,
      \]
      for some $B_{k,l}$, where $N_{I_{k,l} (L_2)} (n)$ denotes the number of words of length $n$ in
      $\Sigma_{I_{k,l} (L_2)}$.
    \end{lemma}

    In order to make use of Lemma \ref{lem:embed}, choose $L_2 > L_3$.

    \subsection{Integrability of the densities}
    The conditional measures of $\mu_\sbr^t$ on unstable manifolds are absolutely continuous
    with respect to Lebesgue measure. We will prove that the conditional measures on the
    stable manifolds are almost surely 
    absolutely continuous with respect to Lebesgue measure. The local product structure of
    $\mu_\sbr^t$ then implies that $\mu_\sbr^t$ is absolutely continuous with respect to 
    Lebesgue measure.

    Take $x \in \Lambda = \pi (\hat{\Lambda})$ and let
    $W_r^{t,\stable} (y,x) = \{z \in W^{t,\stable} (x) \mid d(y,z) \leq r \}$.
    The derivative of $\mu_\sbr^{t,\stable,x}$ at $y$ is the limit
    \[
      D(\mu_\sbr^{t,\stable,x}, y) = \liminf_{r\rightarrow 0} \frac{\mu_\sbr^{t,\stable,x} 
      (W_r^{t,\stable} (y,x))}{2 r}.
    \]
    If the function $D(\mu_\sbr^{t,\stable,x}, y)$ is integrable on $W^{t,\stable} (x)$ then the measure
    $\mu_\sbr^{t,\stable,x}$ is absolutely continuous with respect to Lebesgue measure. 

    Let $k$ be fixed. We want to prove that for \mbox{a.e.} $t \in I_k \cap J_2$
    \begin{equation} \label{Dint}
      \int_{\Omega_t} \int_{\Omega_t} D(\mu_\sbr^{t, \stable,x}
      |_{\Omega_t}, y) \diff{\mu_\sbr^{t,
        \stable,x}} (y) \diff{\mu_\sbr^t} (x) < \infty.
    \end{equation}
    This implies that the measure
    $\mu_\sbr^{t, \stable,x}$ restricted to the set
    $\Omega_t$ is absolutely continuous for
    \mbox{a.e.} $x \in \Omega_t$. Since the conditional measure on the unstable manifolds
    are absolutely continuous with respect to Lebesgue measure, this
    implies that $\mu_\sbr^t |_{\Omega_t}$ is absolutely
    continuous with respect to Lebesgue measure. Since $\mu_\sbr^t
    (\Omega_t) > 0$, ergodicity then implies that this also
    holds for the measure $\mu_\sbr^t$. Since $k$ is arbitrary this implies that
    $\mu_\sbr^t$ is absolutely continuous with respect to Lebesgue for a.e. $t \in I \cap J_2$.

    Fatou's lemma implies that in order to prove \eqref{Dint} it
    suffices to prove that
    \[
      \liminf_{r \rightarrow 0} \frac{1}{r} \int_{\Omega_t} \int_{\Omega_t}
      \mu_\sbr^{t, \stable,x} (\Omega_t \cap
      W_r^{t,\stable} (y,x)) \diff{\mu_\sbr^{t,
        \stable,x}} (y) \diff{\mu_\sbr^t} (x) < \infty.
    \]
    We may rewrite this as
    \begin{equation} \label{intfixedlambda}
      \liminf_{r \rightarrow 0} \frac{1}{r} \int_{\Omega_t} \int_{\Omega_t}
      \int_{\Omega_t} \chi_{\{|y_1-z_1|<r \}}
      \diff{\mu_\sbr^{t, \stable,x}} (z)
      \diff{\mu_\sbr^{t, \stable,x}} (y)
      \diff{\mu_\sbr^t} (x) < \infty.
    \end{equation}
    To prove that this holds for a.e. $t \in I_{k} \cap J_2$ we prove that for any $l$
    \begin{multline} \label{int}
      \liminf_{r \rightarrow 0} \frac{1}{r} \int_{I_{k,l} \cap J_2} \int_{\Omega_t} \int_{\Omega_t}
      \int_{\Omega_t} \chi_{\{|y_1-z_1|<r \}}
      \diff{\mu_\sbr^{t, \stable,x}} (z)
      \diff{\mu_\sbr^{t, \stable,x}} (y)
      \diff{\mu_\sbr^t} (x)
      \diff{t} \\ < \infty.
    \end{multline}
    This then implies that $\mu_\sbr^t
    \ll \nu$ for \mbox{a.e.} $t \in I_{k} \cap J_2$. Instead of
    proving \eqref{int} we use that $\mu_\sbr^t = \hat{\mu}_\sbr^t \circ
    \pi^{-1}$ and
    prove the equivalent condition
    \begin{multline} \label{hatint}
      \liminf_{r \rightarrow 0} \frac{1}{r} \int_{I_{k,l} \cap J_2}
      \int_{\hat{\Omega}_t} \int_{\hat{\Omega}_t} \int_{\hat{\Omega}_t}
      \chi_{\{|y_1-z_1|<r \}} 
      \diff{\hat{\mu}_\sbr^{t, \stable,\hat{x}}} (\hat{z})
      \diff{\hat{\mu}_\sbr^{t, \stable,\hat{x}}} (\hat{y})
      \diff{\hat{\mu}_\sbr^{t}} (\hat{x})
      \diff{t} \\ < \infty.
    \end{multline}

    Recall that $\rho_t : \Sigma_t \rightarrow \hat{\Lambda}_t$ maps sequences 
    in $\Sigma_t$ to points in $\hat{\Lambda}_t$ in the natural way.
    Put $\mu_\Sigma^t = \hat{\mu}_\sbr^t \circ \rho_t$ 
    and rewrite \eqref{hatint} as
    \begin{multline} \label{sigmaint}
      \liminf_{r \rightarrow 0} \frac{1}{r} \int_{I_{k,l} \cap J_2}
      \int_{\hat{\Omega}_t} \int_{\rho_t^{-1}(\hat{\Omega}_t)} \int_{\rho_t^{-1} (\hat{\Omega}_t)}
      \chi_{\{|\rho (\{i_n\}) - \rho (\{j_n\})|_1 < r \}} \\
      \diff{\mu_\Sigma^{t, \stable,\hat{x}}} (\{i_n\})
      \diff{\mu_\Sigma^{t, \stable,\hat{x}}} (\{j_n\})
      \diff{\hat{\mu}_\sbr^{t}} (\hat{x})
      \diff{t}  < \infty,
    \end{multline}
    where $| \cdots |_1$ denotes the difference in the first coordinate.

    Embed all subshifts $\Sigma_{t}$, $t \in I_{k,l} \cap J_1$ into the larger subshift $\Sigma_{I_{k,l}}$
    according to
    Lemma \ref{lem:embed}. The measures $\mu_\Sigma^t$ extend from $\Sigma_t$ to $\Sigma_{I_{k,l}}$
    in a natural way since $\Sigma_t$ is a subset of $\Sigma_{I_{k,l}}$. A cylinder in $\Sigma_{I_{k,l}}$ 
    will be denoted by
    \[
      {}_l [ \{i_n\} ]_m = {}_l [i_{l} \cdots i_m ]_m = \{ \{j_n\} \in \Sigma_{I_{k,l}} \mid j_n = i_n, 
      n=l,\ldots,m \}.
    \]

    To prove \eqref{sigmaint} we estimate the quantity
    \begin{multline} \label{Tr}
      T_r (\{\hat{\Omega}_t \mid t \in I_{k,l} \cap J_2\}) = \\
      \sum_{L>L_2} \hspace{\stretch{1}} \sum_{\substack{{}_{-L} [ i_{-L},\ldots, i_0 ]_0 \\
      \subset \Sigma_{I_{k,l}}}} \hspace{\stretch{1}} \sum_{\substack{1 \leq l_1, l_2 \leq a\\ l_1 \ne l_2}}
      \int_{I_{k,l} \cap J_2} \int_{\hat{\Omega}_t}
      \int_{\rho_t^{-1} (\hat{\Omega}_t) \cap [ l_1, i_{-L}, \ldots, i_0 ] }
      \int_{\rho_t^{-1} (\hat{\Omega}_t) \cap [ l_2, i_{-L}, \ldots, i_0 ] } \\
      \chi_{\{|\rho_t (\{i_n\}) - \rho_t (\{j_n\}) |_1<r \}}
      \diff{\mu_\Sigma^{t, \stable,\hat{x}}} (\{j_n\})
      \diff{\mu_\Sigma^{t, \stable,\hat{x}}} (\{i_n\})
      \diff{\hat{\mu}_\sbr^{t}} (\hat{x})
      \diff{t}
    \end{multline}
    and show that $T_r < \eta r$ for all $r>0$ and some constant $\eta$. This implies \eqref{sigmaint} as follows.
    The product $\Sigma_{I_{k,l}} \times \Sigma_{I_{k,l}}$ can be written as
    \begin{multline*}
      \Sigma_{I_{k,l}} \times \Sigma_{I_{k,l}} = \bigcup_L \ \ \bigcup_{\substack{{}_{-L} [ i_{-L},\ldots, i_0 ]_0 \\
      \subset \Sigma_{I_{k,l}}}} \ \ \bigcup_{\substack{1 \leq l_1, l_2 \leq a\\ l_1 \ne l_2}} \\
      {}_{-L-1} [ l_1, i_{-L}, \ldots, i_0 ]_0 \times {}_{-L-1} [ l_2, i_{-L}, \ldots, i_0 ]_0,
    \end{multline*}
    i.e. $\Sigma_{I_{k,l}} \times \Sigma_{I_{k,l}}$ is the union over $L$ of the set of pair of sequences
    with the first $L$ letters equal. This implies that
   \begin{align*}
      \liminf_{r \rightarrow 0} \frac{1}{r} \int_{I_{k,l} \cap J_2}
      \int_{\hat{\Omega}_t} & \int_{\rho_t^{-1}(\hat{\Omega}_t)} \int_{\rho_t^{-1} (\hat{\Omega}_t)}
      \chi_{\{|\rho (\{i_n\}) - \rho (\{j_n\})|_1 < r \}} \\
      & \hspace{2cm} \diff{\mu_\Sigma^{t, \stable,\hat{x}}} (\{i_n\})
      \diff{\mu_\Sigma^{t, \stable,\hat{x}}} (\{j_n\})
      \diff{\hat{\mu}_\sbr^{t}} (\hat{x})
      \diff{t} \\
      & = \liminf_{r \rightarrow 0} \frac{1}{r} T_r (\{\hat{\Omega}_t \mid t \in I_{k,l} \cap J_2\}),
    \end{align*}
    so \eqref{Tr} implies \eqref{sigmaint}.

    It remains to show \eqref{Tr}. This will be done in \ref{subsec:finalstep}. To do this, the estimate in
    \ref{subsec:power} is needed.

    \subsection{An estimate on power series} \label{subsec:power}
      The expression $|\rho_t (\{i_n\}) - \rho_t (\{j_n\}) |_1$ appearing in \eqref{Tr} can be expressed
      as a power series. The following estimate on this power series is an important part in proving $\eqref{Tr}$.

      If $(x_1,x_2,x_3) \in \hat{\Lambda}_t$ and $\{i_n\} \in \Sigma_{I_{k,l}}$ is the sequence such that
      $\rho_t (\{i_n\}) = (x_1, x_2, x_3)$, then it is easy to see that
      \[
        x_1 = \sum_{n=1}^\infty \prod_{l=1}^{n-1} (t \lambda_{i_{l-n}}) u_{i_{-n}} = 
        \sum_{n=1}^\infty \prod_{l=1}^{n-1} \frac{\lambda_{i_{l-n}}}{\lambda_{\max}} u_{i_{-n}} 
        \bigl( \lambda_{\max} t \bigr) ^{n-1}.
      \]
      So the expression in the brackets of integrand in \eqref{Tr} can be rewritten in the form
      \begin{equation} \label{powerseries}
        \Biggl|  \sum_{n=L+1}^\infty \Bigl(\prod_{l=1}^{n-1} \lambda_{i_{l-n}} u_{i_{-n}}
        - \prod_{l=1}^{n-1} \lambda_{j_{l-n}} u_{j_{-n}} \Bigr) t^{n-1} \Biggr| < r.
      \end{equation}
      This is equivalent to
      \begin{multline} \label{transseries}
        \Biggl|  1
        + \sum_{n=L+2}^\infty \Bigl(\prod_{l=1}^{n-1} \frac{\lambda_{i_{l-n}}}{\lambda_{\max}} u_{i_{-n}}
        - \prod_{l=1}^{n-1} \frac{\lambda_{j_{l-n}}}{\lambda_{\max}} u_{j_{-n}} \Bigr)  
        \frac{(t\lambda_{\max} )^{n-1-L}}{u_{l_1}-u_{l_2}} \Biggr| \\
        < \frac{r}{|u_{l_1}-u_{l_2}|} \frac{1}{\prod_{l=1}^L (t \lambda_{i_{l-n}})},
      \end{multline}
      The coefficients of $(t \lambda_{\max})^{n-1-L}$
      in the sum in \eqref{transseries} are bounded by
      \begin{multline*}
        \max_{\{i_n\}, \{j_n\},n} \Bigl(\prod_{l=1}^{n-1} 
        \frac{\lambda_{i_{l-n}}}{\lambda_{\max}} u_{i_{-n}}
        - \prod_{l=1}^{n-1} \frac{\lambda_{j_{l-n}}}{\lambda_{\max}} u_{j_{-n}} \Bigr) 
        \frac{1}{u_{l_1}-u_{l_2}} \\
        \leq \frac{\max \{ |u_i - u_j|, |u_i| \}}{\min \{ |u_i - u_j| \mid u_i \ne u_j \}} =: C,
      \end{multline*}
      so an application of Corollary \ref{cor:trans} with $C$ as above
      shows that
      \begin{multline} \label{transestimate}
        \Bigl| \Bigl\{ t \lambda_{\max} \mid \Bigl| 
        \sum_{n=L}^\infty (\prod_{l=1}^{n-1} \lambda_{i_{l-n}} u_{i_{-n}}
        - \prod_{l=1}^{n-1} \lambda_{j_{l-n}} u_{j_{-n}}) t^{n-1} \Bigr| < r 
        \Bigr\} \Bigr| \\ < \delta^{-1} \frac{(t_0 \lambda_{\min})^{-L}}{\min \{ |u_i - u_j| \mid
        u_i \ne u_j \}} r
        = E \lambda_{\max} (t_0 \lambda_{\min})^{-L} r,
      \end{multline}
      if $t_1 \lambda_{\max} < \min \Bigl\{\frac{1}{{2C}}, 0.68 \Bigr\}$.

    \subsection{The final step} \label{subsec:finalstep}
      Put $F (t, \{i_n\}, \{j_n\}) = \chi_{\{|\rho_t (\{i_n\}) - \rho_t (\{j_n\}) |_1<r \}}$
      and rewrite \eqref{Tr} as
      \begin{multline} \label{Tr2}
        T_r (\{\hat{\Omega}_t \mid t \in I_{k,l} \cap J_2\}) =  \\
        \sum_{L>L_2} \hspace{\stretch{1}} \sum_{\substack{{}_{-L} [ i_{-L},\ldots, i_0 ]_0 \\
        \subset \Sigma_{I_{k,l}}}} \hspace{\stretch{1}} \sum_{\substack{1 \leq l_1, l_2 \leq a\\ l_1 \ne l_2}}
        \int_{I_{k,l} \cap J_2} \int_{\hat{\Omega}_t}
        \int_{\rho_t (\hat{\Omega}_t) \cap [ l_1, i_{-L}, \ldots, i_0] }
        \int_{\rho_t (\hat{\Omega}_t) \cap [ l_2, i_{-L}, \ldots, i_0] } \\
        F (t, \{i_n \}, \{ j_n \}) 
        \diff{\mu_\Sigma^{t, \stable,\hat{x}}} (\{j_n\})
        \diff{\mu_\Sigma^{t, \stable,\hat{x}}} (\{i_n\})
        \diff{\hat{\mu}_\sbr^{t}} (\hat{x})
        \diff{t}.
      \end{multline}
      To estimate the quantity in \eqref{Tr2} we want to change the order of integration to integrate
      with respect to $t$ first and then use
      the estimate \eqref{transestimate}. This can not be done immediately because the other integrals
      depends on $t$. To get around this problem the function $F (t, \{i_n \}, \{ j_n \})$ will be bounded
      by a function $G (t, \{i_n \}, \{ j_n \})$ which is constant on cylinders. More precisely
      \begin{lemma} \label{lem:G}
        For each pair of cylinders $[ l_1, i_{-L}, \ldots, i_0]$ and $[ l_2, i_{-L}, \ldots, i_0]$
        appearing in \eqref{Tr2} there is a partition into finitely many cylinders
        \[
          [l_1, i_{-L}, \ldots, i_0] = \bigcup_{p=1}^{n_\alpha} S_\alpha^p,\
          [l_2, i_{-L}, \ldots, i_0] = \bigcup_{q=1}^{n_\beta} S_\beta^q
        \]
        and functions $G_{S_\alpha^p, S_\beta^q} (t)$ with
        \[
          G_{S_\alpha^p, S_\beta^q} (t) \geq F (t, \{i_n\}, \{j_n\}), 
          \ \ \forall t \in I_{k,l} \cap J_1,\ \forall \{i_n\} \in S_\alpha^p, \ 
          \forall \{j_n\} \in S_\beta^q
        \]
        and
        \begin{equation} \label{Gint}
          \int_{I_{k,l}} G_{S_\alpha^p, S_\beta^q} (t) \diff{t} < 2 E \lambda_{\max} (t_0 \lambda_{\min})^{-L} r.
        \end{equation}
      \end{lemma}

      \begin{proof}
        For fixed $\hat{y}$ and $\hat{z}$ the estimate \eqref{transestimate} implies that 
        \[
          |\{t \mid F (t, \{i_n\}, \{j_n\}) > 0\}| 
          < E \lambda_{\max} (t_0 \lambda_{\min})^{-L} r
        \]
        and
        \begin{equation} \label{rlength}
          |\{t \mid |\rho_t (\{i_n\}) - \rho_t (\{j_n\}) |_1<2r \}| < 
          2 E \lambda_{\max} (t_0 \lambda_{\min})^{-L} r.
        \end{equation}
        The function $(t, \{i_n\}, \{j_n\}) \mapsto |\rho_t (\{i_n\}) - \rho_t (\{j_n\}) |_1$
        depends continuously
        on $t$, $\{i_n\}$ and $\{j_n\}$. Choose $L'$ so large that
        \[
          \frac{2 \max\{u_i - u_j\} (t_1 \lambda_{\max})^{L'+3}}{1 - t_1 \lambda_{\max}} < r.
        \]
        The sets $[l_1, i_{-L}, \ldots, i_0]$ and $[l_1, i_{-L}, \ldots, i_0]$
        can be partitioned into cylinders of length $L' + L + 3$
        \begin{align*}
          [l_1, i_{-L}, \ldots, i_0] &= \bigcup_{p=1}^{n_\alpha} S_\alpha^p = 
          \bigcup_{p=1}^{n_\alpha} [\alpha_{-L'} (p), \ldots, \alpha_{0} (p), l_1, i_{-L}, \ldots, i_0],\\
          [l_2, i_{-L}, \ldots, i_0] &= \bigcup_{q=1}^{n_\beta} S_\beta^q = 
          \bigcup_{q=1}^{n_\beta} [\beta_{-L'} (q), \ldots, \beta_{0} (q), l_2, i_{-L}, \ldots, i_0].
        \end{align*}
        Note that $n_\alpha$ and $n_\beta$ can be bounded uniformly by $n_\alpha, n_\beta \leq a^{L'+1}$.

        As in \eqref{powerseries} the expression
        $||\rho_t(\{i_n^{(1)}\}) - \rho_t (\{j_n^{(1)}\})|_1 - | \rho_t(\{i_n^{(2)}\}) - \rho_t (\{j_n^{(2)}\})|_1|$
        can be written as
        \begin{multline*}
          ||\rho_t(\{i_n^{(1)}\}) - \rho_t (\{j_n^{(1)}\})|_1 - | \rho_t(\{i_n^{(2)}\}) -
          \rho_t (\{j_n^{(2)}\})|_1| \\
          = \Biggl| \Biggl| \sum_{n=1}^\infty \Bigl(\prod_{l=1}^{n-1} \lambda_{i_{l-n}^{(1)}} u_{i_{-n}^{(1)}}
          - \prod_{l=1}^{n-1} \lambda_{j_{l-n}^{(1)}} u_{j_{-n}^{(1)}} \Bigr) t^{n-1} \Biggr|\\
          - \Biggl| \sum_{n=1}^\infty \Bigl(\prod_{l=1}^{n-1} \lambda_{i_{l-n}^{(2)}} u_{i_{-n}^{(2)}}
          - \prod_{l=1}^{n-1} \lambda_{j_{l-n}^{(2)}} u_{j_{-n}^{(2)}} \Bigr) t^{n-1} \Biggr| \Biggr|
        \end{multline*}
        If $\{i_n^{(1)}\}, \{i_n^{(2)}\} \in S_\alpha^p$, $\{j_n^{(1)}\}, \{j_n^{(2)}\} \in S_\beta^q$
        for some $p$ and $q$, then since $\{i_n^{(1)}\}$ and $\{i_n^{(2)}\}$ respectively $\{j_n^{(1)}\}$ and 
        $\{j_n^{(2)}\}$ are equal on the first $L' + L + 3$ letters, the first $L' + L + 3$ terms in this powerseries
        are zero. Hence
        \begin{align} \label{cylinderestimate}
          ||\rho_t(\{i_n^{(1)}\}) &- \rho_t (\{j_n^{(1)}\})|_1 - | \rho_t(\{i_n^{(2)}\}) -
          \rho_t (\{j_n^{(2)}\})|_1| = \nonumber \\
          &= \Biggl| \Biggl| \sum_{n=L+L'+4}^\infty \Bigl(\prod_{l=1}^{n-1} \lambda_{i_{l-n}^{(1)}} u_{i_{-n}^{(1)}}
            - \prod_{l=1}^{n-1} \lambda_{j_{l-n}^{(1)}} u_{j_{-n}^{(1)}} \Bigr) t^{n-1} \Biggr| \nonumber \\
          &\phantom{= \Biggl|} - \Biggl| \sum_{n=L+L'+4}^\infty \Bigl(\prod_{l=1}^{n-1} 
            \lambda_{i_{l-n}^{(2)}} u_{i_{-n}^{(2)}}
            - \prod_{l=1}^{n-1} \lambda_{j_{l-n}^{(2)}} u_{j_{-n}^{(2)}} \Bigr) t^{n-1} \Biggr| \Biggr| \nonumber \\
          &\leq \sum_{n=L+L'+4}^\infty 2 \max \{u_i - u_j\} (\lambda_{\max} t)^{n-1} \nonumber \\
          &\leq \frac{2 \max \{u_i - u_j\} (\lambda_{\max} t_1)^{L+L'+3}}{1 - \lambda_{\max} t_1} < r.
        \end{align}
        The last inequality follows from the choice of $L'$ that
        $\frac{2 \max\{u_i - u_j\} (t_1 \lambda_{\max})^{L'+3}}{1 - t_1 \lambda_{\max}} < r$.

        For each pair $S_\alpha^p$ and $S_\beta^q$, take $\{i_n^{(p)}\} \in S_\alpha^p$ and
        $\{j_n^{(q)}\} \in S_\beta^q$. Put
        \[
          G_{S_\alpha^p, S_\beta^q} (t)= \chi_{\{|\rho_t (\{i_n^{(p)}\}) - \rho_t (\{j_n^{(q)}\}) |_1< 2r \}}.
        \]
        Then the estimates \eqref{rlength} and \eqref{cylinderestimate} imply that
        \[
          G_{S_\alpha^p, S_\beta^q} (t) \geq F (t, \{i_n\}, \{j_n\}), 
          \ \forall t \in I_{k,l} \cap J_1, \forall \{i_n\} \in S_\alpha^p,
          \forall \{j_n\} \in S_\beta^q
        \]
        and
        \[
          \int_{I_{k,l}} G_{S_\alpha^p, S_\beta^q} (t) \diff{t} < 2 E \lambda_{\max} (t_0 \lambda_{\min})^{-L} r.
          \qedhere
        \]
      \end{proof}

      Lemma \ref{lem:G} is used to estimate the integrals in \eqref{Tr2} in the following way.
      \begin{align*}
        \mathcal{I} :&= \int_{I_{k,l} \cap J_2} \int_{\hat{\Omega}_t}
        \int_{\rho_t^{-1} (\hat{\Omega}_t) \cap [l_1, i_{-L}, \ldots, i_0] }
        \int_{\rho_t^{-1} (\hat{\Omega}_t) \cap [l_2, i_{-L}, \ldots, i_0] } \\
        &\hspace{2cm} F_{L,i_{-L}, \cdots, i_0} (t, \{i_n\}, \{j_n\}) 
        \diff{\mu_\Sigma^{t, \stable,\hat{x}}} (\{j_n\})
        \diff{\mu_\Sigma^{t, \stable,\hat{x}}} (\{i_n\})
        \diff{\hat{\mu}_\sbr^{t}} (\hat{x})
        \diff{t} \\
        &= \sum_{p,q} \int_{I_{k,l} \cap J_2} \int_{\hat{\Omega}_t}
        \int_{\rho_t^{-1} (\hat{\Omega}_t) \cap S_\alpha^p}
        \int_{\rho_t^{-1} (\hat{\Omega}_t) \cap S_\beta^q } \\
        &\hspace{1cm} F_{L,i_{-L}, \cdots, i_0} (t, \{i_n\}, \{j_n\}) 
        \diff{\mu_\Sigma^{t, \stable,\hat{x}}} (\{j_n\})
        \diff{\mu_\Sigma^{t, \stable,\hat{x}}} (\{i_n\})
        \diff{\hat{\mu}_\sbr^{t}} (\hat{x})
        \diff{t} \\
        &\leq \sum_{p,q} \int_{I_{k,l} \cap J_2} \int_{\hat{\Omega}_t}
        \int_{\rho_t^{-1} (\hat{\Omega}_t) \cap S_\alpha^p}
        \int_{\rho_t^{-1} (\hat{\Omega}_t) \cap S_\beta^q } \\
        &\hspace{2cm} G_{S_\alpha^p,S_\beta^q} (t) 
        \diff{\mu_\Sigma^{t, \stable,\hat{x}}} (\{j_n\})
        \diff{\mu_\Sigma^{t, \stable,\hat{x}}} (\{i_n\})
        \diff{\hat{\mu}_\sbr^{t}} (\hat{x})
        \diff{t} \\
        & = \sum_{p,q} \int_{I_{k,l} \cap J_2} \int_{\hat{\Omega}_t} G_{S_\alpha^p,S_\beta^q} (t) \\
        &\hspace{2cm}  
        \mu_\Sigma^{t, \stable,\hat{x}} 
        (\rho_t^{-1} (\hat{\Omega}_t) \cap S_\alpha^p)
        \mu_\Sigma^{t, \stable,\hat{x}} (\rho_t^{-1} (\hat{\Omega}_t)
        \cap S_\beta^q) 
        \diff{\hat{\mu}_\sbr^{t}} (\hat{x}) \diff{t}.
      \end{align*}
      This can be estimated by
      \begin{multline*}
        \mathcal{I} \leq \max_{p,q} \biggl( \int_{I_{k,l} \cap J_1} G_{S_\alpha^p,S_\beta^q} (t) \diff{t} \biggr) \\
        \sup_t \biggl( \sum_{p,q} \int_{\hat{\Omega}_t}
        \mu_\Sigma^{t, \stable,\hat{x}} 
        (\rho_t^{-1} (\hat{\Omega}_t) \cap S_\alpha^p)
        \mu_\Sigma^{t, \stable,\hat{x}} (\rho_t^{-1} (\hat{\Omega}_t) 
        \cap S_\beta^q) \diff{\hat{\mu}_\sbr^{t}} (\hat{x}) \biggr).
      \end{multline*}
      The maximum is estimated by \eqref{Gint} and the sum can be eliminated using that $\{S_\alpha^p\}$
      and $\{S_\beta^q\}$ are partitions of the cylinders $[l_1, i_{-L}, \ldots, i_0]$ and
      $[l_1, i_{-L}, \ldots, i_0]$ respectively:
      \begin{multline*}
        \mathcal{I} \leq 2 E \lambda_{\max} (t_0 \lambda_{\min})^{-L} r
        \sup_t \biggl( \int_{\hat{\Omega}_t} 
        \mu_\Sigma^{t, \stable,\hat{x}} 
        (\rho_t^{-1} (\hat{\Omega}_t) \cap [l_1, i_{-L}, \ldots, i_0]) \\
        \mu_\Sigma^{t, \stable,\hat{x}} (\rho_t^{-1} (\hat{\Omega}_t)
        \cap [l_2, i_{-L}, \ldots, i_0])
        \diff{\hat{\mu}_\sbr^{t}} (\hat{x}) \biggr).
      \end{multline*}
      The measure of the cylinders is estimated with \eqref{smb}. 
      \begin{align*}
        \mathcal{I} &\leq 2 E \lambda_{\max} (t_0 \lambda_{\min})^{-L} r 
        \sup_t \biggl( \int_{\hat{\Omega}_t} 
        A_k^2 (\gamma_{\min} - \varepsilon_2)^{-2 (L+1)} \diff{\hat{\mu}_\sbr^{t}} (\hat{x}) \biggr) \\
        &\leq 2 E \lambda_{\max} (t_0 \lambda_{\min})^{-L} r 
        A_k^2 (\gamma_{\min} - \varepsilon_2)^{-2 (L+1)}.
      \end{align*}
      Thus
      \begin{multline*}
        T_r (\{\hat{\Omega}_t \mid t \in I_{k,l} \cap J_2\}) \\ \leq \sum_{L = L_2}^\infty
        \sum_{\seq{i}} \sum_{l=1}^m 
        \sum_{\substack{1 \leq l_1, l_2 \leq a \\ l_1\ne l_2}} 2 E r (t_0 \lambda_{\min})^{-L}
        A_k^2 (\gamma_{\min} - \varepsilon_2)^{-2 (L+1)} \\
        \leq \sum_{L=L_2}^\infty 2 E r B_{k,l} (\gamma_{\max} + \varepsilon_2 + \varepsilon_4)^L m 
        (a^2-a) (t_0 \lambda_{\min})^{-L} 
        A_k^2 (\gamma_{\min} - \varepsilon_2)^{-2L} \\
        = \sum_{L=L_2}^\infty 2 E A_k^2 B_{k,l} m (a^2 - a) r
        \biggl( \frac{\gamma_{\max} + \varepsilon_2 + \varepsilon_4}
        {t_0 \lambda_{\min} (\gamma_{\min} - \varepsilon_2)^2} \biggr)^L.
      \end{multline*}
      If $\frac{t_0 \lambda_{\min} \gamma_{\min}^2}{\gamma_{\max}} > 1$ then it is possible to choose 
      $\varepsilon_2$ and $\varepsilon_4$ so small that
      \[
        \tau = \frac{\gamma_{\max} + \varepsilon_2 + \varepsilon_4}{t_0 \lambda_{\min}
        (\gamma_{\min} - \varepsilon_2)^2} < 1.
      \]
      Then $T_r \leq 2 E A_k^2 B_{k,l} m (a^2 - a) \tau^{L_2} \frac{1}{1 - \tau} r$. This implies that 
      $\mu_\sbr^t$ is absolutely continuous with respect to Lebesgue measure for a.e. 
      $t \in I \cap J_2$. 
      Let $\varepsilon_0, \varepsilon_1, \varepsilon_3 \rightarrow 0$. 
      Then $\nu (I \cap J_2) > \nu(I) - \varepsilon_0 - \varepsilon_1 - \varepsilon_3 \rightarrow \nu (I)$
      and this shows that
      $\hat{\mu}_\sbr^t$ is absolutely continuous with
      respect to Lebesgue measure for a.e. $t \in I_k$.

  \section{Application to fat Belykh maps}
    The Belykh maps are defined as follows.
    Let $-1 < k < 1$ and put
    $K = [-1, 1]^2$, $K_1 = \{(x_1,x_2) \in (-1,1)^2 \mid  x_2 > k x_1\}$ and
    $K_2 = \{(x_1,x_2)\in (-1,1)^2 \mid x_2 < k x_1\}$. The discontinuity set is
    $N = \{(x_1,x_2) \mid x_2 = k x_1\}$ and the border is 
    $M = \{(x_1,x_2) \in K \mid x_1 = \pm 1 \text{ or } x_2 = \pm 1\}$.
    The Belykh maps $f_{\lambda, \gamma, k} : K_1 \cup K_2 \rightarrow K$ are defined by
    \begin{align*}
      &f_{\lambda, \gamma, k} |_{K_1} (x_1,x_2) = (\lambda x_1 + (1 - \lambda),
      \gamma x_2 - (\gamma -1)), \\
      &f_{\lambda, \gamma, k} |_{K_2} (x_1,x_2) = (\lambda x_1 - (1 - \lambda),
      \gamma x_2 + (\gamma -1)),
    \end{align*}
    where $0 < \lambda < 1$ and $1 < \gamma \leq \frac{2}{1 + |k|}$.

    These maps were first introduced for $\lambda < \frac{1}{2}$ by Belykh in \cite{Belykh} as a model 
    of a Poincar\'{e} map from phase synchronisation. It was investigated in the case
    $\lambda < \frac{1}{2}$
    in \cite{Pesin} and \cite{Sataev}. Schmeling and Troubetzkoy studied in
    \cite{Schmeling-Troubetzkoy} the 
    Belykh maps when $\lambda > \frac{1}{2}$ and called this case the fat Belykh map.

    The Belykh maps satisfy the conditions (A1) and (A2). Here
    \[
      C = 1 \text{ and } C < f_4 (x) \text{ if } x < Q_4 = 0.61
    \]
    and we conclude:
    \begin{theorem}
      Let $P = \{ (\lambda, \gamma, k) \mid \gamma \lambda > 1,\ \lambda < 0.61\}$.
      For Lebesgue almost all  $(\gamma, \lambda, k) \in P$ the fat Belykh map
      $f_{\lambda, \gamma, k}$ has an absolutely 
      continuous invariant measure.
    \end{theorem}

  \section{Decay of correlations}
    Applying Young's scheme from \cite{Young} Chernov proved in \cite{Chernov} the exponential decay of correlations
    for H\"{o}lder continuous functions for a class of piecewise hyperbolic systems with singularities in arbitrary
    dimensions. This result can be used in the following way.

    Let $H_\eta = \{\phi: K \rightarrow \R \mid \exists C : |\phi(x) - \phi(y)| \leq C d(x,y)^\eta, \
    \forall x,y \in K \}$
    be the set of H\"{o}lder continuous functions on $K$.
    \begin{theorem} \label{the:correlations}
      Assume that $f: K \rightarrow K$ satisfies the assumptions (A1) and (A2) and assume that
      $(f^n, \mu_\sbr)$ is ergodic for every $n \geq 1$.
      For every $\eta > 0$ there exists a constant $\theta \in (0,1)$ such that for every $\phi, \psi \in H_\eta$
      there is a constant $C(\phi, \psi)$ such that
      \[
        \biggl| \int_K \phi \circ f^n \cdot \psi \diff{\mu_\sbr} -
        \int_K \phi \diff{\mu_\sbr} \int_K \psi \diff{\mu_\sbr} \biggr| \leq C \theta^n
      \]
      for all $n \in \N$.
    \end{theorem}

    \begin{proof}
      Lift $f : K \rightarrow K$ to $\hat{f} : \hat{K} \rightarrow \hat{K}$ as in Section \ref{subsec:notation}.
      The lift $\hat{f}$ is just a model for the natural extension of $f$. The natural extension is
      ergodic if and only if the original system is ergodic, see Theorem 1 in Chapter 10, \S 4 of
      \cite{Cornfeld-Fomin-Sinai}. Hence $(\hat{f}^n, \hat{\mu}_\sbr)$ is ergodic for every $n \geq 1$.

      The function $\phi : K \rightarrow \R$ is lifted to $\hat{\phi} : \hat{K} \rightarrow \R$ by
      $\hat{\phi} = \phi \circ \pi$ and $\hat{\psi}$ is lifted in the same way.

      Note that $\hat{\phi}, \hat{\psi} \in \hat{H}_\eta =\{ \hat{\phi} : \hat{K} \rightarrow \R \mid
      \exists C : |\hat{\phi} (x) - \hat{\phi} (y)| \leq C d(x,y)^\eta, \ \forall x,y \in \hat{K} \}$ and
      \begin{align*}
        & \int_K \phi \circ f^n \cdot \psi \diff{\mu_\sbr} = \int_{\hat{K}} \hat{\phi} \circ \hat{f}^n
        \cdot \hat{\psi} \diff{\hat{\mu}_\sbr}, \\
        & \int_K \phi \diff{\mu_\sbr} = \int_{\hat{K}} \hat{\phi} \diff{\hat{\mu}_\sbr}, \
        \int_K \psi \diff{\mu_\sbr} = \int_{\hat{K}} \hat{\psi} \diff{\hat{\mu}_\sbr}.
      \end{align*}
      Theorem 1.1 in \cite{Chernov} states that $(\hat{f}, \hat{\mu}_\sbr)$ has exponential decay of correlations
      so this implies that $(f, \mu_\sbr)$ has exponential decay of correlations.
    \end{proof}

\end{document}